\documentclass[11pt]{article}

\usepackage{latexsym,amssymb}
\usepackage{amsmath}
\usepackage{amscd}
\usepackage{hyperref}

\input xypic
\xyoption{all}

\newcommand{\rar}{\rightarrow}
\newcommand{\lar}{\longrightarrow}

\newcommand{\llar}{-\kern-5pt-\kern-5pt\longrightarrow}
\newcommand{\lllar}{-\kern-5pt-\kern-5pt\llar}

\newtheorem{Theorem}{Theorem}[section]
\newtheorem{Lemma}[Theorem]{Lemma}
\newtheorem{Corollary}[Theorem]{Corollary}
\newtheorem{Proposition}[Theorem]{Proposition}

\newtheorem{Conjecture}[Theorem]{Conjecture}

\newtheorem{Remark}[Theorem]{Remark}
\newtheorem{Example}[Theorem]{Example}
\newtheorem{Definition}[Theorem]{Definition}
\newtheorem{Question}[Theorem]{Question}

\def\coker{\mbox{\rm coker}}
\def\demo{\noindent{\bf Proof. }}
\def\depth{\mbox{\rm depth }}
\def\ds{\displaystyle}
\def\edeg{\mbox{\rm edeg}}
\def\Ext{\mbox{\rm Ext}}
\def\Sym{\mbox{\rm Sym}}
\def\gr{\mbox{\rm gr}}
\def\reg{\mbox{\rm reg}}
\def\red{\mbox{\rm red}}

\def\height{\mbox{\rm height }}
\def\Hom{\mbox{\rm Hom}}

\def\length{{\lambda}}

\def\m{\mathfrak{m}}

\def\QED{\hfill$\Box$}

\def\Rees{\mbox{$\mathcal{R}$}}

\def\Symi{\mbox{\rm Sym}}

\def\AA{{\mathbf A}}

\def\BB{{\mathbf B}}
\def\CC{{\mathbf C}}

\def\RR{{\mathbf R}}
\def\SS{{\mathbf S}}

\def\LL{{\mathbf L}}
\def\xx{{\mathbf x}}
\def\TT{{\mathbf T}}
\def\ff{{\mathbf f}}

\def\g2{{\mathbf g}}
\def\hh{{\mathbf h}}
\def\xx{{\mathbf x}}

\def\ttt{{\mathbf t}}
\def\H{{\mathrm H}}

\def\m{{\mathfrak m}}
\def\p{{\mathfrak p}}

\begin{document}

\title{{\sc Extremal Rees Algebras\\ \date{\today}}
\vspace{-1mm}}\footnotetext[1]{2010 AMS {\it Mathematics Subject
Classification}: 13B10, 13D02, 13H10,  13H15, 14E05. }
\footnotetext[2]{The second author is partially supported by CNPq, Brazil. The last author is partially supported by the NSF. }
\footnotetext[3]{Corresponding author: Wolmer V. Vasconcelos.}

\author{
{\normalsize\sc Jooyoun Hong}
\vspace{-0.75mm}\\
{\small Department of Mathematics}\vspace{-1.4mm} \\
{\small Southern Connecticut State University}\vspace{-1.4mm}\\
{\small 501 Crescent Street}\vspace{-1.4mm}\\
{\small New Haven, CT 06515-1533, U.S.A.}\vspace{-1.4mm}\\
{\small e-mail: {\tt hongj2@southernct.edu}}\vspace{4mm}
\and
{\normalsize\sc Aron Simis}
\vspace{-0.75mm}\\
{\small Departamento de Matem\'atica}\vspace{-1.4mm} \\
{\small Universidade Federal de Pernambuco}\vspace{-1.4mm}\\
{\small 50740-560 Recife, PE, Brazil}\vspace{-1.4mm}\\
{\small e-mail: {\tt aron@dmat.ufpe.br}}\vspace{4mm}
\and
{\normalsize\sc Wolmer V. Vasconcelos}
\vspace{-0.75mm}\\
{\small Department of Mathematics}\vspace{-1.4mm} \\
{\small Rutgers University}\vspace{-1.4mm}\\
{\small 110 Frelinghuysen Road}\vspace{-1.4mm}\\
{\small Piscataway, NJ 08854-8019, U.S.A.}\vspace{-1.4mm} \\
{\small e-mail: {\tt vasconce@math.rutgers.edu}}\vspace{4mm}
}

\maketitle

\begin{center}

{\em Dedicated to J\"urgen Herzog for his numerous contributions to Commutative Algebra on the occasion of his 70th birthday} 
\end{center}



\begin{abstract}
\noindent
We study almost complete intersections ideals whose Rees algebras are extremal in the sense that some of their fundamental metrics---depth or relation type---have maximal or minimal values in the class.
The focus is on those 
ideals that lead to {\em almost} Cohen--Macaulay algebras and our treatment is wholly concentrated on 
 the nonlinear relations of the algebras. 
Several classes of such algebras are presented, some  of a combinatorial origin. We offer a different prism to look at 
these questions with accompanying techniques.
The main results are  effective methods to calculate the invariants of these algebras.

\end{abstract}


\noindent {\small {\bf  Key Words and Phrases:} Almost complete intersection, almost Cohen--Macaulay algebra,
  birational mapping, Castelnuovo regularity, extremal Rees algebra, Hilbert function, module of nonlinear relations, Rees algebra, relation type, Sally module.
}

\newpage

\tableofcontents

\section{Introduction}

\noindent
Our goal  is the study of the defining equations of the Rees algebras $\RR[It]$ of classes of almost complete intersection ideals when one of its important metrics--especially depth or reduction number--attains an extreme value in
 the class.
 We are going to show that such algebras occur frequently and develop novel means to identify them.
 As a consequence interesting  properties  of such algebras have been discovered.  
 We argue
 that several questions, while often placed in the general context of Rees algebra theory, may
 be viewed as subproblems in this more narrowly defined environment.
 
  \medskip

Let $\RR$ be a Cohen--Macaulay local ring of dimension $d$, or a polynomial ring 
 $\RR=k[x_1, \ldots, x_d]$ for
$k$ a field.
By an {\em almost complete intersection} we mean an ideal
  $I=(a_1, \ldots, a_g, a_{g+1})$ of codimension $g$ where the subideal $J=(a_1, \ldots, a_g)$ is a 
  complete intersection    and  $a_{g+1}\notin J$.
    By the {\em equations} of $I$ it is meant a free presentation of
 the Rees algebra $\RR[It]$ of $I$,
\begin{eqnarray}\label{presRees1}
 0 \rar \LL \lar \BB = \RR[\TT_1, \ldots, \TT_{g+1}] \stackrel{\psi}{\lar}
\RR[It] \rar 0,  \quad \TT_i \mapsto f_it .
\end{eqnarray}
More precisely,  $\LL$ is the defining ideal of the Rees
algebra of $I$ but we refer to it simply as the {\em ideal of
equations} of $I$. We are particularly interested in establishing the properties of $\LL$
when $\RR[It]$ is Cohen--Macaulay or has almost maximal depth.
This broader view requires a change of focus from $\LL$ to one of its quotients.
 We are going to study some classes of ideals whose 
Rees algebras have these properties. They tend to occur in classes where the reduction number $\red_J(I)$ attains an extremal value.

\medskip

We first set up the framework to deal with properties of $\LL$ by a standard decomposition. We keep the notation of above, $I = (J,a)$. The presentation ideal $\LL$ of $\RR[It]$
 is a graded ideal
$\LL = L_1 + L_2 + \cdots$, where $L_1$ are linear forms in the $\TT_i$ defined by a	matrix $\phi$ of the syzygies of $I$,
$ L_1 = [\TT_1, \ldots, \TT_{g+1}] \cdot \phi $. 
Our basic prism is given by the exact sequence
\[ 0 \rar \LL/(L_1) \lar \BB/(L_1) \lar \RR[It] \rar 0.\] 
Here
$\BB/(L_1)$ is a presentation of the symmetric algebra of $I$ and  
$\SS = \Symi(I)$ is a Cohen--Macaulay ring under very broad conditions, 
including when $I$ is an ideal of finite colength. The emphasis here will be entirely on $T = \LL/(L_1)$, which we call the
{\em module of nonlinear} relations of $I$. The usefulness arises because  of the  fact  exhibited in
the exact sequence
\begin{eqnarray}\label{presRees2}
0 \rar T \lar \SS \lar \RR[It] \rar 0.
\end{eqnarray}
\begin{itemize}
\item[{$\bullet$}] [Proposition~\ref{canoseq}] $T$ is a Cohen--Macaulay $\SS$--module if and only if $\depth
\RR[It]\geq d$.
\end{itemize}
Rees algebras with this property will be called {\em almost Cohen--Macaulay} (aCM for short). We note that $\LL$ carries very  
different
kind of information than $T$ does. The advantage lies in the flexibility of treating Cohen--Macaulay modules 
 Cohen--Macaulay ideals: the means to test for Cohen-Macaulayness in modules are more plentiful than in ideals.
  An elementary example lies in the proof of:

\begin{itemize}
\item[{$\bullet$}] [Theorem~\ref{reducedsymi}]    Suppose  that $\RR$ is a Cohen--Macaulay  local ring and  $I$ is
an $\m$--primary 
 almost  complete intersection such that  $\SS=\Sym(I)$ is reduced.  Then $\RR[It]$ is almost Cohen--Macaulay.
 \end{itemize}

\medskip

We shall now discuss our more technical results. Throughout $(\RR,\m)$ is a Cohen--Macaulay local ring of dimension 
$d$ (where we include rings of polynomials and the ideals are homogeneous), and $I$ is an almost complete intersection $I = (J,a)$ of finite colength. 

\medskip

One technique we bring in to the treatment of the equations of Rees algebras is the theory of the Sally module. It gives
a very direct relationship between the Cohen--Macaulayness of $T$ and of the Sally module $S_J(I)$ of $I$ relative 
to $J$. $S_J(I)$ gives also a quick connection between the Castelnuovo regularity and the relation type of
$\RR[It]$ and those of $S_J(I)$. A criterion of Huckaba (\cite{Huc96}) (of which we give a quick proof for
completeness) gives a method to test the Cohen--Macaulayness of $T$ in terms of the values of the first
Hilbert coefficient $e_1(I)$ of $I$ (Theorem~\ref{Huckaba}).
 It is particularly well-suited for the case when $I$ is generated by
homogeneous polynomials defining a birational mapping for then the value of $e_1(I)$ is known. 

\medskip


%

Our approach to the estimation for $\nu(T)$, the minimum number of generators of $T$, passes through the determination
of an effective formula for  $\deg \SS$, the multiplicity of $\SS$:

\begin{itemize}
\item[{$\bullet$}] [Theorem~\ref{degSymibis}] If $I$ is generated by forms of degree $n$, then
\[ \deg \SS = \sum_{j=0}^{d-1} n^j + \lambda(I/J).\]
\end{itemize}
The summation accounts for  $\deg \RR[It]$, according to \cite{HTU}, so $\deg T = \lambda(I/J)$. This is a number that will control the number of generators of $T$, and therefore of $\LL$, whenever $\RR[It]$ is almost Cohen--Macaulay.
It achieves 
 the goal  to find estimates for the number of generators of $\LL$ and of its {\em relation type}, that is
\[ \mbox{\rm reltype}(I) = \inf\{n \mid \LL = (L_1, L_2, \ldots, L_n)\}.\]

Two other metrics of interest, widely studied 
 for homogeneous ideals but not limited to them, are the following. 
One   
   seeks to bound the {\em saturation exponent} of $\LL/(L_1)$ (which was 
 introduced in \cite{syl2} and has a simple ring-theoretic explanation as the index of nilpotency of $\SS$),
\[ \mbox{\rm sdeg}(I) = \inf\{s \mid \m^s \LL \subset (L_1)\},\]
and the other
  is the degree of the special fiber $\mathcal{F}(I)$ of $\RR[It]$, also called the {\em elimination degree} of $I$,
\[ \edeg(I) = \deg \mathcal{F}(I)= \inf\{s \mid L_s \not \subset \m\BB\}.\]

While $\mbox{\rm retype}(I)$ is the most critical of these numbers, the other two are significant because they are often
found linked to the syzygies of $I$. Our notion of extremality will cover the supremum or infimum values of these degrees in a given class of ideals but also their relationship to the cohomology of $\RR[It]$ as expressed by the depth of the algebra.

\medskip

 Two classes of almost Cohen--Macaulay algebras arise from certain homogeneous ideals. First, we show that 
binary ideals with one linear syzygy have this property. This has been proved by several authors. We offer a very
short proof using the technology of the Sally module (Proposition~\ref{aCMofbin}). It runs for a few lines and gives no details of the projective resolution of that algebra besides the fact that it has the appropriate length. We include it because we have found no similar technique in the literature.
The proof
structure, a simple combinatorial obstruction to the aCM property, is used repeatedly to examine the occurrence of
the property amongst ideals generated by quadrics in $4$--space.

\medskip

A different class of algebras are those associated to monomials. These ideals have the form $I=(x^{\alpha}, y^{\beta}, 
z^{\gamma}, x^ay^bz^c)$. We showed that
\begin{itemize}
\item[{$\bullet$}] [Proposition~\ref{nnn111}] The ideals $(x^n, y^n, z^n, xyz)$, $n\geq 3$, and $(x^n, y^n, z^n, w^n, xyzw)$, $n\geq 4$, 
have almost Cohen--Macaulay Rees algebras.
\end{itemize}

We expect these statements are still valid in higher dimensions. Our proofs were computer-assisted as we used
Macaulay2 (\cite{Macaulay2}) to derive deeper heuristics.

\section{Approximation complexes and almost complete intersections}

A main source of extremal Rees algebras lie in the construction of approximation complexes. We quickly recall them and some of their main properties.

\subsection{The $\mathcal{Z}$--complex}

 These are complexes derived from Koszul
complexes, and arise as follows (for details, see \cite{HSV1}, \cite{HSV83},
 \cite[Chapter 4]{alt}). Let $\RR$ be a commutative ring, $F$ a free $\RR$-module of rank $n$,
with a basis  $\{e_1, \ldots, e_n\}$,
 and
$\varphi: F\rar \RR$ a homomorphism. The exterior algebra $\bigwedge F$ of
$F$ can be endowed with a differential
\[ \partial_{\varphi} : \bigwedge^{r} F
{\lar} \bigwedge^{r-1}F,
\]
\[ \partial_{\varphi}(v_1\wedge v_2\wedge \cdots \wedge v_{r}) =
\sum_{i=1}^{r} (-1)^i\varphi(v_i) (v_1\wedge \cdots \wedge
\widehat{v_i} \wedge  \cdots \wedge v_{r}).  \]
The complex $\mathbb{K}(\varphi)=\{ \bigwedge F, \partial_{\varphi}\}$ is called the {\em
Koszul complex} of $\varphi$. \index{Koszul complex} Another notation
for it is: Let $\xx=\{\varphi(e_1), \ldots, \varphi(e_n)\}$, denote
the Koszul complex by $\mathbb{K}(\xx)$.

\medskip

Let $\SS = S(F)= \Sym(F) = \RR[\TT_1, \ldots, \TT_n]$, and consider the
exterior algebra of $F\otimes_{\RR} S(F)$. It can be viewed as a Koszul
complex obtained from $\{\bigwedge F, \partial_{\varphi}\}$ by change
of scalars $\RR \rar \SS$, and another complex defined by the
$\SS$-homomorphism
\[ \psi: F \otimes_{\RR}S(F) \lar S(F), \quad \psi(e_i)= \TT_i.
\]
The two differentials $\partial_{\varphi}$ and $\partial_{\psi}$
satisfy
\[ \partial_{\varphi}\partial_{\psi}+
\partial_{\psi}\partial_{\varphi}=0,\]
which leads directly to the construction of several
complexes.

\begin{Definition}{\rm  Let $\mathbf{Z}$, $\mathbf{B}$ and $\mathbf{H}$ be the modules of
cycles, boundaries  and
 the homology of $\mathbb{K}(\varphi)$.
\begin{itemize}
\item The $\mathcal{Z}$-complex of $\varphi$ is
$ \mathcal{Z}=\{ \mathbf{Z}\otimes_{\RR} \SS, \partial\}$ 
\[ 0 \rar Z_n\otimes \SS[-n] \lar \cdots \lar Z_1\otimes \SS[-1] \lar \SS \rar 0,\]
where $\partial$ is the differential induced by $\partial_{\phi}$.
\item The $\mathcal{B}$-complex of $\varphi$ is
the sub complex of $ \mathcal{Z}$ 
\[ 0 \rar B_n\otimes \SS[-n] \lar \cdots \lar B_1\otimes \SS[-1] \lar \SS \rar 0.\]
\item The $\mathcal{M}$-complex of $\varphi$ is
$\mathcal{M}=\{ \mathbf{H}\otimes_{\RR} S, \partial\}$
\[ 0 \rar \H_n\otimes
 \SS[-n] \lar \cdots \lar \H_1\otimes \SS[-1] \lar \H_0 \otimes \SS \rar 0,\]
where $\partial$ is the differential induced by $\partial_{\psi}$.
\end{itemize}
}\end{Definition}

These are complexes of graded modules over the polynomial ring
$\SS=\RR[\TT_1, \ldots, \TT_n]$.

\begin{Proposition} Let $I = \varphi(F)$. Then
\begin{enumerate}
\item[{\rm (i)}]  The homology of $\mathcal{Z}$  and of $\mathcal{M}$ depend only on $I$;
\item[{\rm (ii)}] $ \H_0(\mathcal{Z})= \Sym(I)$.
\end{enumerate}
\end{Proposition}

%
%
%


\subsubsection*{Acyclicity}
The homology of the Koszul complex $\mathbf{K}(\varphi)$ is not fully
independent of $I$, for instance, it depends on the number of
generators. An interest here is the ideals whose $\mathcal{Z}$ complexes are acyclic.
%
%

We recall a broad setting that gives rise to almost complete intersections with Cohen--Macaulay symmetric algebras. For 
a systematic examination of the notions here we refer to \cite{HSV83}. It is centered on one approximation complex associated to an ideal, the so-called $\mathcal{Z}$--complex.
An significant interest for us is the following.
\begin{Theorem}[{\cite[Theorem 10.1]{HSV83}}] Let $\RR$ be a Cohen-Macaulay local ring and let $I$ be an
ideal of positive height. Assume:
\begin{itemize}
\item[{\rm (a)}] $\nu(I_{\mathfrak{p}})\leq \height \mathfrak{p}+1$ for
every \ $\mathfrak{p} \supset I$;

\item[{\rm (b)}] $\depth (\H_i)_{\mathfrak{p}} \geq \height
\mathfrak{p}-\nu(I_{\mathfrak{p}})+1 $ for every \ $\mathfrak{p}\supset
I$ and every $0\leq i \leq \nu(I_{\mathfrak{p}})-\height
I_{\mathfrak{p}} $.
\end{itemize}
Then
\begin{itemize}
\item[{\rm (i)}] The complex $\mathcal{Z}$ is acyclic.
\item[{\rm (ii)}] $\Sym(I)$ is a  Cohen-Macaulay ring.
\end{itemize}
\end{Theorem}

\begin{Corollary}\label{Zaci} Let $\RR$ be a Cohen-Macaulay local ring of dimension $d \geq 1$ and let $I$ be an
almost complete intersection.
The complex $\mathcal{Z}$ is acyclic and $\Sym(I)$ is a
Cohen-Macaulay algebra in the following cases:
\begin{itemize} 
\item[{\rm (i)}]  \mbox{\rm (See also \cite{Rossi83})} $I$ is $\m$-primary. In this case $\Sym(I)$ has Cohen--Macaulay type $d - 1$. 
\item[{\rm (ii)}]  $\height I= d-1$. 
Furthermore if $I$ is
generically a complete intersection then $I$ is of linear type.
\item[{\rm (iii)}]  $\height I= d-2$ and
 $\depth \RR/I\geq 1$. Furthermore if $\nu(I_{\mathfrak{p}})\leq
  \height \mathfrak{p}$ for $I\subset \mathfrak{p}$ then $I$ is of
  linear type.
\end{itemize}
\end{Corollary}

A different class of ideals with Cohen--Macaulay symmetric algebras is treated in \cite{Johnson}.

\subsection{The canonical presentation}

Let $\RR$ be a Cohen--Macaulay local domain of dimension $d\geq 1$ and let $I$ be an almost complete intersection as in Corollary~\ref{Zaci}. The ideal of equations $\LL$ can be studied  in two stages: $(L_1)$ and $\LL/(L_1)= T$:
\begin{eqnarray} \label{canoseque} 
 0 \rar T \lar \SS = \BB/(L_1)= \Sym(I) \lar \RR[It] \rar 0.\end{eqnarray}
We will argue that this exact sequence is very useful. Note that $\Sym(I)$ and $\RR[It]$ have dimension $d+1$, and that $T$ is the $\RR$--torsion submodule of $\SS$.
Let us give some of its properties.

\begin{Proposition} \label{canoseq} Let $I$ be an ideal as above.
\begin{enumerate}

\item[{\rm (i)}]  $(L_1)$ is a Cohen-Macaulay ideal of $\BB$.

\item[{\rm (ii)}] $T$ is a Cohen--Macaulay $\SS$--module if and only if $\depth \RR[It]\geq d$.

\item[{\rm (iii)}] If $I$ is $\m$--primary then $\mathcal{N}=T\cap \m\SS$ is the nil radical of 
$\SS$ and $\mathcal{N}^s = 0$ if and only if $\m^sT =0$. This is equivalent to saying that $\mbox{\rm sdeg}(I)$ is the index of nilpotency of $\Symi(I)$.

\item[{\rm (iv)}] $T=\mathcal{N} + \mathcal{F}$, where $\mathcal{F}$ is a lift in $\SS$ of the relations in $\SS/\m \SS$ of the special
fiber ring $\mathcal{F}(I)=\RR[It] \otimes \RR/\m$.  In particular if $\mathcal{F}(I)$ is a hyersurface ring,
$T = (f, \mathcal{N})$.

\end{enumerate}

\end{Proposition}

\demo (i) comes from Corollary~\ref{Zaci}.
\medskip

\noindent (ii)  In the defining sequence of $T$,
\[ 0 \rar (L_1) \lar \LL \lar T \rar 0,\]
since $(L_1)$ is a Cohen--Macaulay ideal of codimension $g$, as an $\BB$--module, we have $\depth (L_1) = d+2$, while
$\depth \LL = 1 + \depth \RR[It]$. It follows that $\depth T =\min \{d+1,  1 +\depth \RR[It]\}$.
Since $T$ is a module of Krull dimension $d+1$ it is a Cohen--Macaulay module if and only if 
$\depth \RR[It] \geq d$.

\medskip

\noindent (iii) $\m\SS$ and $T$ are both minimal primes and for large $n$, $\m^nT=0$. Thus $T$ and $\m\SS$ are the only minimal 
primes of $\SS$, $\mathcal{N} = \m \SS \cap T$. 
To argue the equality of the two indices of nilpotency, let $n$ be such that $\m^nT=0$. The ideal $\m^n\SS + T$
has positive codimension, so contains regular elements since $\SS$ is Cohen--Macaulay. Therefore to show
 \[\m^s T=0 \Longleftrightarrow \mathcal{N}^s =0\]
it is enough to multiply both expressions by $\m^n\SS + T$. The verification is immediate.

\medskip

\noindent (iv) Tensoring the sequence (\ref{canoseque}) by $\RR/\m$ gives the exact sequence
\[ 0 \rar \m\SS\cap T/\m T = \mathcal{N}/\m T \lar T/\m T \lar \SS/\m \SS \lar \mathcal{F}(I) \rar 0. \]
By Nakayama Lemma we may ignore $\m T$ and recover $T$ as asserted. \QED

\bigskip

The main intuition derived from these basic observations is (ii): Whenever methods are developed to study the equations of
$\RR[It]$ when this algebra is Cohen--Macaulay, should apply [hopefully] in case they are almost Cohen--Macaulay.

\begin{Remark}{\rm If $I$ is not $\m$--primary
but still satisfies one of the other conditions of Corollary~\ref{Zaci},
the nilradical $\mathcal{N}$ of $\SS$ is given by $T\cap N_0\SS$, where $N_0$ is the intersection of the minimal
primes $\p$ for which $I_{\p}$ is not of linear type.

}
\end{Remark}

\subsection{Reduced symmetric algebras}

\begin{Proposition}\label{linkofmax} Let $\RR$ be a Gorenstein local domain of dimension $d$ and let $J$ be a parameter
ideal. If $J$ contains two minimal generators in $\m^2$,
then the Rees algebra of $I=J\colon \m$ is Cohen--Macaulay and 
$\LL = ( L_1, \ff)$ for some quadratic form $\ff$. 
\end{Proposition}

\demo
The equality $I^2 = JI$  comes from \cite{CPV1}. The Cohen--Macaulayess of $\RR[It]$ is a general argument 
(in \cite{CPV1} and probably elsewhere). Let $\ff$ denote the quadratic form 
\[ \ff = \TT_{d+1} ^2 + \mbox{\rm lower terms}.\]

Let us show that $\m T=0$.
Reduction modulo $\ff$ can be used to present any element in $\LL$ as 
\[ F =  \TT_{d+1} \cdot A + B\in \LL, \]
where $A$ and $B$ are forms in $\TT_1,\ldots, \TT_d$. Since $I = J:\m$,
 any element in $\m \TT_{d+1}$ is equivalent, modulo 
$L_1$, to a linear form in the other variables.
Consequently  
\[\m F \subset (L_1, \RR[\TT_1, \ldots, \TT_d]) \cap \LL \subset (L_1),\]
as desired.

By Proposition~\ref{canoseq}, we have the exact sequence
\[ 0\rar T \lar \SS/\m \SS \lar \mathcal{F}(I) \rar 0.\]
But $T$ is a maximal Cohen--Macaulay $\SS$--module, and so it is also a maximal Cohen--Macaulay $\SS/\m \SS$--module as
well. It follows that $T$ is generated by a monic polynomial that divides the image of $\ff$ in $\SS/\m\SS$. 
It is now clear that $T=(\ff)\SS$. \QED


\begin{Corollary} For the ideals above $\Sym(I)$ is reduced.
\end{Corollary}


We now discuss a generalization, but since we are still developing the examples, we are somewhat informal.

\begin{Corollary}\label{sdegs}
Suppose the syzygies of $I$ are contained in $\m^s\BB$ and that $\m^s\LL \subset  (L_1)$.  We have the exact sequence
\begin{eqnarray} \label{canoseque2}
 0 \rar T \lar \SS/\m^s\SS \lar \RR[It]\otimes \RR/\m^s = \mathcal{F}_s(I) \rar 0.
\end{eqnarray}
If $\RR[It]$ is almost Cohen--Macaulay, $T$ is a Cohen--Macaulay module that is an ideal of the polynomial
ring $\CC = \RR/(\m^s)[\TT_1, \ldots, \TT_{d+1}]$, a ring of multiplicity
${s+d-1}\choose{d}$.
Therefore we have that $\nu(T) \leq {{s+d-1}\choose {d}}$.
\end{Corollary}

Note that also here $\mathcal{F}_s(I)$ is Cohen--Macaulay. We wonder whether $\mathcal{F}(I)$ is Cohen--Macaulay.

\bigskip

Let $(\RR, \m)$ be a Cohen--Macaulay local ring and $I$ an almost complete intersection as in Corollary~\ref{Zaci}.
We examine the following surprising fact.

\begin{Theorem} \label{reducedsymi}    Suppose  that $\RR$ is a Cohen--Macaulay  local ring and  $I$ is
an $\m$--primary 
 almost  complete intersection such that  $\SS=\Sym(I)$ is reduced.  Then 
 $\RR[It]$ is an almost Cohen-Macaulay algebra.
 \end{Theorem}
 
 \demo 
  Since 
$0 =\mathcal{N} = T\cap \m \SS$,   on one hand from (\ref{canoseque}) we have that $T$ satisfies the condition $S_2$ of Serre, that is
\[ \depth T_{P} \geq \inf\{2, \dim T_P\}\]
for every prime ideal $P$ of $\SS$. On the other hand, 
  from (\ref{canoseque2}) $T$ is an ideal of the polynomial ring $\SS/\m \SS$. It follows that
$T= (\ff )\SS$, and consequently $\depth \
 \RR[It] \geq d$.
\QED

\begin{Example}{\rm  If $\RR = \mathbb{Q}[x,y]/(y^4-x^3)$, $J = (x)$ and $I = J: (x,y)= (x, y^3)$, $\depth \RR[It] = 1$.
}\end{Example}

There are a number of immediate observations.

\begin{Corollary} If $I$ is an ideal as in {\rm Theorem~\ref{reducedsymi}}, then the special fiber ring $\mathcal{F}(I)$ is Cohen-Macaulay.
\end{Corollary}

\begin{Remark}{\rm  If $I$ is an almost complete intersection as in (\ref{Zaci}) and its radical is a regular prime
ideal $P$, that is  $\RR/P$ is regular local ring, the same assertions will apply if $\Symi(I)$ is reduced.
}\end{Remark}

\section{Almost Cohen--Macaulay algebras}

We begin our treatment of the properties of an ideal when its
 Rees algebra $\AA =\RR [It]$ is
almost Cohen-Macaulay. We first describe a large class of examples.

\subsection{Direct links of Gorenstein ideals}

We briefly outline a broad class of extremal Rees algebras.
Let $(\RR, \m)$ be a 
Gorenstein local ring of dimension $d\geq 1$. A natural source of almost complete intersections in $\RR$ direct links of Gorenstein ideals. That 
is, let $K$ be a Gorenstein ideal of $\RR$ of codimension $s$, that is $\RR/K$ is a Gorenstein ring of dimension $d-s$.
If $J = (a_1, \ldots, a_s)\subset K$ is a complete intersection of codimension $s$, $J\neq K$, $I = J:K$ is an almost complete intersection, $I = (J, a)$. Depending on $K$, sometimes these ideals come endowed with very good
properties. Let us recall one of them.

\begin{Proposition}\label{sourceofacis} Let $(\RR, \m)$ be a Noetherian local ring of dimension $d$.
\begin{enumerate}
\item[{\rm (i)}] {\rm (\cite[Theorem 2.1]{CPV1})} Suppose $\RR$ is a Cohen--Macaulay local ring and  let $\p$ be a prime ideal of codimension $s$ such that $\RR_{\p}$ is a Gorenstein ring and let $J$ be a complete intersection of 
codimension $s$ contained in $\p$. Then for $I= J:\p$ we have $I^2= JI$ in the following two cases: {\rm (a)} $\RR_{\p}$ is not a regular local ring; {\rm (b)} if $\RR_{\p}$ is a regular local ring two of the elements $a_i$ belong to $\p^{(2)}$. 

\item[{\rm (ii)}] {\rm (\cite[Theorem 3.7]{CHV})}  	Suppose $J$ is an irreducible $\m$--primary ideal. Then {\rm (a)} either there exists a minimal set of generators $\{x_1, \ldots, x_d\}$ of $\m$ such that $J = (x_1,\ldots, x_{d-1}, {x_d}^r)$, or {\rm (b)}
$I^2 = JI$ for $I = J:\m$.

\end{enumerate}
\end{Proposition}

\medskip

The following criterion is a global version of Corollary~\ref{F2Cor}

\begin{Proposition}\label{rednumone} Let $\RR$ be a Gorenstein local ring and $I=(J,a)$ an almost complete intersection {\rm[}when we write $I=(J,a)$ we always mean that $J$ is a reduction{\rm]}. If $I$ is an unmixed ideal {\rm[}height unmixed{\rm]} then
$\red_J(I) \leq 1$ if and only if $J:a = I_1(\phi)$.
\end{Proposition}

\demo Since the ideal $JI$ is also unmixed, to check the equality $J:a= I_1(\phi)$ we only need to check at the minimal primes of $I$ (or, of $J$, as  they are the same). Now Corollary~\ref{F2Cor} applies. \QED

\bigskip

If in Proposition~\ref{sourceofacis} $\RR$ is a Gorenstein local ring and $I$ is a Cohen--Macaulay ideal, their
associated graded rings are Cohen--Macaulay, while the 
Rees algebras are also Cohen--Macaulay if $\dim \RR\geq 2$.

\begin{Theorem} \label{reesoflink} Let $\RR$ be a Gorenstein local ring and $I$ a Cohen--Macaulay ideal that is an  almost complete intersection. If $\red_J(I)\leq 1$ then in the canonical representation
\[ 0 \rar T \lar \SS \lar \RR[It] \rar 0,\]
\begin{itemize}
\item[{\rm (i)}] If $\dim \RR \geq 2$ $\RR[It]$ is Cohen--Macaulay.

\item[{\rm (ii)}] $T$ is a Cohen--Macaulay module over $\SS/(I_1(\phi))\SS$, in particular
\[ \nu(T) \leq \deg \RR/I_1(\phi).\]
\end{itemize}
\end{Theorem}

\begin{Example}{\rm 
Let $\RR=k[x_1, \ldots, x_d]$, $k$ an algebraically closed field, and let $\p$ be a homogeneous prime ideal of codimension $d-1$. Suppose $J= (a_1, \ldots, a_{d-1})$ 
is a complete intersection of codimension $d-1$ with at least two generators in $\p^2$. Since $\RR/\p$ is regular,
$I=J:\p$ is an almost complete intersection and $I^2 = JI$. Since $\p$ is a complete intersection, say
$\p = (x_1-c_1x_{d} , x_2-c_2x_d, \ldots, x_{d-1} -c_{d-1}x_d)$, $c_i\in k$,  we write  the matrix equation $J= \AA \cdot \p$, where
$\AA$ is a square matrix of size $d-1$. This is the setting where the Northcott ideals occur, and therefore
$I = (J, \det \AA)$.

\medskip

By Theorem~\ref{reesoflink}(ii), $\nu(T) \leq \deg (\RR/\p) = 1$. Thus $\LL$ is generated by the syzygies of $I$ (which are well-understood) plus a quadratic equation.
}\end{Example}

\subsection{Metrics of aCM Rees algebras}

Let $(\RR, \m)$ be a Cohen--Macaulay local ring of dimension $d$ and let $I$ be an almost complete intersection of 
finite colength. We assume that $I=(J,a)$, where $J$ is a minimal reduction of $I$. These assumptions will hold 
for the remainder of the section. We emphasize that they apply to the case when $\RR$ is a polynomial ring over
a field and $I$ is a homogeneous ideal.

\medskip

In the next statement we highlight the information about the equations of $I$ that is a direct consequences of the aCM hypothesis. In the next segments we begin to obtain the required data in an explicit form. As for
 notation, $\BB = \RR[\TT_1, \ldots, \TT_{d+1}]$ and 
  for a graded 
 $\BB$-module $A$, $\deg(A)$ denotes the  multiplicity relative to the maximal homogeneous ideal $\mathcal{M}$ of $\BB$,
 $\deg(A)= \deg( \gr_{\mathcal{M}}(A))$. In actual computations $\mathcal{M}$ can be replaced by a reduction. For instance, if $E$ is a graded $\RR$--module and $A=E\otimes_{\RR} \BB$, picking a reduction  $J$ for $\m$ gives the reduction 
 $\mathcal{N} = (J, \TT_1, \ldots, \TT_{d+1})$ of $\mathcal{M}$. It will follow that $\deg(A) = \deg(E)$.

\begin{Theorem} \label{aCM1} If the algebra $\RR[It]$ is almost Cohen--Macaulay, in the canonical sequence
\[ 0 \rar T \lar \SS \lar \RR[It] \rar 0\]
\begin{itemize}
\item[{\rm (i)}] $\reg(\RR[It]) = \red_J(I) + 1$.

\item[{\rm (ii)}] $\nu(T) \leq \deg(\SS) - \deg(\RR[It])$.

\end{itemize}
\end{Theorem}

\demo (i) follows from Corollary~\ref{Sallyrel}. As for (ii), since $T$ is a Cohen--Macaulay module, $\nu(T) \leq \deg(T)$. \QED

\bigskip

The goal is 
 to find $\deg(T)$, $\deg(\RR[It])$ and $\deg(\SS)$ in terms of more direct metrics of $I$. This will be answered in Theorem~\ref{degSymi}.

\subsection*{Cohen--Macaulayness of the Sally module.}

Fortunately there is a simple criterion to test whether $\RR[It]$ is an aCM algebra: It is so
if and only if  it satisfies the Huckaba Test:
\[ e_1(I) = \sum_{j\geq 1}\length(I^j/JI^{j-1}). \]
Needless to say, this is exceedingly effective if you already know
$e_1(I)$,  in particular there is no need to determine the equations of $\RR[It]$ for the
purpose.

\medskip

Let $\RR$ be a Noetherian ring, $I$ an ideal and $J$ a reduction of $I$. The Sally module of
$I$ relative to $J$, $S_J(I)$, is defined by the exact sequence of $\RR[Jt]$--modules
\[ 0 \rar I \RR[Jt] \lar I \RR[It] \lar S_J(I) = \bigoplus_{j\geq 2} I^j/IJ^{j-1} \rar 0.\]
  The definition applies more broadly to other filtrations. We refer the reader to \cite[p. 101]{icbook} for
  a discussion. Of course this module depends on the chosen reduction 
 $J$, but its Hilbert function and its depth are independent of $J$.
There are extensions of this construction to more general reductions--and we employ one below. 
 
 \medskip

If $\RR$ is a Cohen--Macaulay local ring and $I$ is $\m$--primary with a minimal reduction, $S_J(I)$ plays
a role in mediating among properties of $\RR[It]$.

 \begin{Proposition} \label{Sallyelem} Suppose $\RR$ is a Cohen--Macaulay local ring of dimension $d$. Then
 \begin{enumerate}
 \item[{\rm (i)}] If $S_J(I) = 0$ then $\gr_I(\RR)$ is Cohen-Macaulay.
 \item[{\rm (ii)}] If $S_J(I)\neq 0$ then $\dim S_J(I) = d$.
 \end{enumerate}
 \end{Proposition}

Some of the key properties of the Sally module  are in
display in the next result (\cite[Theorem 3.1]{Huc96}). It converts the property of $\RR[It]$ being almost Cohen--Macaulay into the property of $S_J(I)$ being Cohen--Macaulay.

\begin{Theorem}[Huckaba Theorem] \label{Huckaba}
 Let $(\RR,\m)$ be a Cohen--Macaulay local ring of
dimension $d \geq 1$ and $J$ a
parameter ideal. Let $\mathcal{A}=\{I_n, n\geq 0\}$ be an filtration
of $\m$-primary ideals such that $J\subset I_1$ and $\BB=\RR[I_nt^n,
n\geq 1]$ is $\AA=\RR[Jt]$-finite. 
Define the Sally module $S_{\BB/\AA}$ 
of $\BB$ relative to $\AA$ by the exact sequence
\[ 0 \rar I_1 \AA \lar I_1\BB \lar S_{\BB/\AA}\rar 0.\]
Suppose $S_{\BB/\AA}\neq 0$.
Then
\begin{enumerate}
\item[{\rm (i)}]  $e_0(S_{\BB/\AA})=e_1(\BB)-\lambda(I_1/J) \leq \sum_{j\geq 2}\lambda(I_j/JI_{j-1})$.
\item[{\rm (ii)}]  The following conditions are equivalent:
\begin{itemize}
\item[{\rm (a)}] $S_{\BB/\AA}$ is Cohen-Macaulay;

\item[{\rm (b)}] $\depth \gr_{\mathcal{A}}(\RR)\geq d-1$;

\item[{\rm (c)}] $e_1(\BB)= \sum_{j \geq 1}\lambda(I_j/JI_{j-1})$;

\item[{\rm (d)}] $\RR[It]$ is almost Cohen--Macaulay.

\end{itemize}

\end{enumerate}

\end{Theorem}

\demo If $J=(\xx)=(x_1,\ldots, x_d)$,
$S_{\BB/\AA}$ is a finite module over the ring $\RR[\TT_1, \ldots,
\TT_d]$, $\TT_i \rar x_it$. Note that
\[\lambda( S_{\BB/\AA}/\xx S_{\BB/\AA})=
\sum_{j\geq 2}\lambda(I_j/JI_{j-1}),
\]
which shows the first assertion.

\medskip

For the equivalencies, first note that equality means that the first
Euler characteristic $\chi_1(\xx;S_{\BB/\AA})$ vanishes, which
by Serre's theorem (\cite[Theorem 4.7.10]{BH}) says that $S_{\BB/\AA}$ is Cohen--Macaulay.
The final assertion comes from the formula for the multiplicity of
$S_{\BB/\AA}$ in terms of $e_1(\BB)$ (\cite[Theorem 2.5]{icbook}).
\QED


\subsection*{Castelnuovo regularity}

The Sally module encodes also  information about the Castelnuovo regularity $\mbox{\rm reg}(\RR[It])$ of
the Rees algebra. The following Proposition and its Corollary are extracted from the literature (\cite{Huck87}, \cite{Trung98}), or proved directly by adding the exact sequence that defines $S_J(I)$ (note that $I\RR[Jt]$ is a maximal Cohen--Macaulay module) to the canonical  sequences relating $\RR[It])$ to $\gr_I(\RR)$ and $\RR$ via $I\RR[It]$ (see
\cite[Section 3]{Trung98}).

\begin{Proposition} \label{Sallyreg} Let $\RR$ be a Cohen--Macaulay local ring, $I$ an $\m$--primary ideal and $J$ a minimal reduction.
Then
\[ \mbox{\rm reg}(\RR[It]) = \mbox{\rm reg}(S_J(I)). \]
In particular
\[ \mbox{\rm reltype}(I) \leq  \mbox{\rm reg}(S_J(I)).\]
\end{Proposition}

\begin{Corollary} \label{Sallyrel} 
If $I$ is an almost complete intersection and $\RR[It]$ is almost Cohen--Macaulay, then
\[ \mbox{\rm reltype}(I) = \red_J(I) + 1.\]
\end{Corollary}

\subsection*{The Sally fiber of an ideal}

To help analyze the problem, we single out an extra structure.
Let $(\RR, \m)$ be a Cohen--Macaulay local ring of dimension $d>0$, $I$ an $\m$--primary ideal and $J$
one of its minimal reductions.

\begin{Definition}[Sally fiber] The Sally fiber of $I$ is the graded module
\[ F(I) = \bigoplus_{j\geq 1} I^j/JI^{j-1}. \]  
\end{Definition}

$F(I)$ is an Artinian $\RR[Jt]$--module whose last non-vanishing component is $I^r/JI^r$, $r=\red_J(I)$. The 
equality $e_1(I) = \length(F(I))$ is the condition for the almost Cohen--Macaulayness of $\RR[It]$. 
We note that $F(I)$ is the fiber of $S_J(I)$ extended by the term $I/J$. 
To
obtain additional control over $F(I)$ we are going to endow it with additional structures in
cases of interest.

\bigskip

Suppose $\RR$ is a Gorenstein local ring, $I=(J,a)$. The modules $F_j = I^j/JI^{j-1}$ are cyclic modules
over the Artinian Gorenstein ring $ \AA = \RR/J:a$. We turn $F(I)$ into a graded module over the  polynomial ring $\AA[s]$ by defining
\[ a^j \in F_j \mapsto  s\cdot a^j = a^{j+1} \in F_{j+1}.\]  
This is clearly well-defined and has $s^r\cdot F(I)=0$. Several of the properties of the $F_n$'s arise from this
representation, for instance the length of $F_j$ are non-increasing. 
Thus $F(I)$ is a graded module over the Artinian Gorenstein
ring $\BB= \AA[s, s^r=0]$.

\medskip

\begin{Remark}\label{Fvasneweq}{\rm 
The variation of the values of  the $F_j$ is connected to the degrees of the generators of $\LL$.  For convenience we set $I=(J,a)$ and $\BB = \RR[u, \TT_1, \ldots, \TT_d]$, with $u$ corresponding to $a$. For example:

\medskip

\begin{itemize}

\item[{\rm (i)}] Suppose that for some $s$, $\ff_s = \lambda(F_s) = 1$. This means that we have $d$ equations of the form
\[\hh_i= x_i u^s + \g2_i \in \LL_s\]
where  $\g2_i\in (\TT_1, \ldots, \TT_{d})\BB_{s-1}$. Eliminating the $x_i$, we derive a nonvanishing monic equation in $\LL$ of degree $d\cdot s$. Thus $\red_J(I) \leq ds -1$.

\medskip

\item[{\rm (ii)}] A more delicate observation,  is that whenever $\ff_s > \ff_{s+1}$ then
there are {\bf fresh} equations in  $\LL_{s+1}$. Let us explain why this happens: 
$\ff_s = \lambda(JI^{s-1}: I^s)$, that is the ideal $L_s$ contains elements of the form
\[ c\cdot u^s + \mathbf{g}, \quad c\in JI^{s-1}:I^s,  \quad  \mathbf{g} \in (\TT_1, \ldots, \TT_d)\BB_{s-1}.\]
Since $\ff_{s+1} < \ff_s$, $JI^s: I^{s+1}$ contains properly $JI^{s-1}:I^s$, which means that we must have elements in
$L_{s+1}$ 
\[ d\cdot u^{s+1} + \mathbf{g}, \quad\]
with $d \notin JI^{s-1}: I^s$ and $\mathbf{g} \in (\TT_1, \ldots, \TT_d)\BB_{s}$. Such elements cannot belong to
$L_s\cdot \BB_1$, so they are fresh generators.

The converse also holds.

\end{itemize}

}\end{Remark}

\subsubsection*{A toolbox}

We first give a simplified version of \cite[Proposition 2.2]{CPV1}.
Suppose $\RR$ is a Gorenstein local ring of dimension $d$. Consider the two
exact sequences.

\[0 \rar J/JI= (\RR/I)^d \lar \RR/JI \lar \RR/J\rar 0\]
 and the  syzygetic sequence
\[ 0 \rar \delta(I) \lar H_1(I) \lar (\RR/I)^{d+1} \lar I/I^2 \rar 0.
\]
The first
gives
\[\length(\RR/JI)= d\cdot \length(\RR/I) +\length(\RR/J),\]
 the other
\[\length(\RR/I^2) =(d+2)\length(\RR/I)-\length(\H_1(I)) + \length
(\delta(I)).\]
Thus
\[\length(I^2/JI)= \length(I/J) -\length(\delta(I))\]
since $\H_1(I)$ is
the canonical module of $\RR/I$. Taking into account the syzygetic
formula in \cite{syl2} we finally have:

\begin{Proposition} \label{F2} Let $(\RR, \m)$ be a
Gorenstein local ring of dimension $d>0$, $J= (a_1, \ldots, a_d)$ a
parameter ideal and $I=(J,a)$ and $a\in \m$. Then
\begin{eqnarray*} \length(I^2/JI) & = &\length(I/J)
-\length(\RR/I_1(\phi))\\ &=& \length(\RR/J:a)  -\length(\RR/I_1(\phi))\\
& = & \length(\RR/J:a) - \length(\Hom(\RR/I_1(\phi), \RR/J:a) \\
& = & \length(\RR/J:a) - \length((J:a): I_1(\phi))/J:a) \\
& = & \length (\RR/(J:a):I_1(\phi)).
\end{eqnarray*} 

\end{Proposition}

Note that in dualizing $\RR/I_1(\phi)$ we made use of the fact that $\RR/J:a$ is a Gorenstein ring.

\begin{Corollary} \label{F2Cor} $I^2= JI$ if and only if $J:a= I_1(\phi)$. In this
case, if $d>1$ the algebra $\RR[It]$ is Cohen--Macaulay.
\end{Corollary}

\begin{Corollary}\label{F3Cor} If $\RR[It]$ is an aCM algebra and $\red_J(I) = 2$, then
$e_1(I) = 2\cdot \lambda(I/J) - \lambda(\RR/I_1(\phi))$.
\end{Corollary}

\begin{Remark}{\rm We could enhance these observations considerably if formulas for $\lambda(JI^2:I^3)$ were to be developed. More precisely, how do the syzygies of $I$ affect $JI^2:I^3$?
}
\end{Remark}

\subsection{Multiplicities and number of relations}

To benefit from Theorem~\ref{aCM1}, we need to have effective formulas for $\deg(\SS)$ and $\deg(\RR[It])$. 
We are going to develop them now.

\begin{Proposition}\label{multirees} Let $\RR=k[x_1, \ldots, x_d]$ and $I$ an almost complete intersection as above,
$I=(f_1, \ldots, f_d, f_{d+1})=(J, f_{d+1})$ generated by 
forms
 of degree $n$. Then
$ \deg(\RR[It])  = \sum_{j=0}^{d-1} n^j.$
\end{Proposition}  

\demo  After an elementary observation, we make use of one of the beautiful multiplicity formulas of \cite{HTU}.
Set $A=\RR[It]$, $A_0 = \RR[Jt]$, $\mathcal{M} = (\m, It)A$ and $\mathcal{M}_0 = (\m, Jt)A_0$. 
Then 
\[ \deg(\gr_{\mathcal{M}_0}(A_0))= \deg(\gr_{\mathcal{M}_0}(A))= \deg(\gr_{\mathcal{M}}(A)), \]
the first equality because $A_0 \rar A$ is a finite rational extension, the second is because $(\m, Jt)A$ is a reduction of $(\m, It)A$. Now we use \cite[Corollary 1.5]{HTU} that gives $\deg(A_0)$. 
\QED

\subsubsection*{The multiplicity of the symmetric algebra}

We shall now prove one of our main results, a formula for $\deg S(I)$ for ideals generated by forms of the same degree.
Let $\RR=k[x_1, \ldots, x_d]$, $I = (\ff) = (f_1, \ldots, f_d, f_{d+1})$ an almost complete intersection generated by forms of degree $n$. At some point we assume, harmlessdly,  that $J = (f_1, \ldots,  f_d)$ is a complete intersection.
There will be a slight change of notation in the rest of this section. We set $\BB = \RR[\TT_1, \ldots, \TT_{d+1}]$ and 
$\SS = \Symi(I)$.

\begin{Theorem}[{\bf Degree Formula}] \label{degSymi} $\deg \SS = \sum_{j=0}^{d} n^j - \lambda(\RR/I)$.
\end{Theorem}

\demo 
Let $\mathbb{K}(\ff) = \bigwedge^{d+1} \RR^{d+1}(-n)$ be the Koszul complex associated to $\ff$,
\[ 0 \rar {K}_{d+1} \rar K_{d} \rar \cdots \rar K_2 \rar K_1 \rar K_0 \rar 0,\]
and consider the associated $\mathcal{Z}$--complex
\[ 0 \rar Z_d\otimes \BB(-d) \rar Z_{d-1} \otimes \BB(-d+1)
\rar \cdots \rar Z_2 \otimes \BB(-2) \rar Z_1 \otimes \BB(-1) \stackrel{\psi}{\rar} \BB \rar 0.\]

This complex is acyclic with $\H_0(\mathcal{Z}) = \SS = \Symi(I)$.  
Now we introduce another complex obtained by replacing $Z_1 \otimes \BB(-1) $ by $B_1 \otimes \BB(-1)$, where
$B_1$ is the module of $1$--boundaries of $\mathbb{K}(\ff)$, 
 followed by
the restriction of $\psi$ to $B_1 \otimes \BB(-1)$.

 \medskip
 
 This defines another acyclic complex, $\mathcal{Z}^*$, actually the $\mathcal{B}$--complex of $\ff$,  and we set $\H_0(\mathcal{Z}^*) = \SS^*$. The relationship
 between $\SS$ and $\SS^*$ is given in the following observation:
 
 \begin{Lemma} $\deg \SS^* = \deg \SS + \lambda(\RR/I)$.
 \end{Lemma}
 
 \demo Consider the natural mapping between $\mathcal{Z} $ and $\mathcal{Z}^*$:
 
 \[
 \diagram
 0 \rto  & Z_{d}\otimes \BB(-d) \rto \dto_{\phi_d} & \cdots \rto  & Z_2\otimes \BB(-2) \rto \dto_{\phi_2} & B_1 \otimes \BB(-1) \rto \dto  & \BB \rto \dto & \SS^* \rto \dto & 0 \\
 0 \rto & Z_{d}\otimes \BB(-d) \rto & \cdots \rto & Z_2\otimes \BB(-2) \rto & Z_1 \otimes \BB(-1) \rto & \BB \rto &
 \SS \rto & 0.
 \enddiagram 
   \]
 The maps $\phi_2, \ldots, \phi_2 $ are isomorphisms while the other maps are defined above. They induce the short exact sequence of modules of dimension $d+1$,
 \[ 0 \rar (Z_1/B_1)\otimes \BB(-1) \lar \SS^* \lar \SS \rar 0.\]
 
Note that $Z_1/B_1= \H_1(\mathbb{K}(\ff))$ is the canonical module of $\RR/I$, and therefore has the same length
as   
$\RR/I$. Finally, by the additivity formula for the multiplicities (\cite[Lemma 13.2]{Eisenbudbook}),
\[ \deg \SS^* = \deg \SS + \lambda(Z_1/B_1),\]
as desired. \QED

\bigskip

We are now give our main calculation of multiplicities.

\begin{Lemma} $\deg \SS^* = \sum_{j=0}^{d} n^j.$

\end{Lemma}

\demo We note that the $\mathcal{Z}^*$--complex is homogeneous for the total degree [as required for the computation of multiplicities] provided the $Z_i$'s and $B_1$ have the same degree. 
We can conveniently write $B_i$ for $Z_i$, $i\geq 2$.
This is clearly the case since they are 
 generated in degree $n$. This is not the case for $Z_1$. However when  $\ff$ is a regular sequence, all the $Z_i$ are equigenerated, an observation we shall make  use of below.

\medskip

Since the modules of $\mathcal{Z}^*$ are homogeneous 
we have that the Hilbert series of $\SS^*$ is given as
\[ H_{\SS^*}(\ttt) = {\frac{\sum_{i=0}^d (-1)^{i} h_{B_i}(\ttt) \ttt^i}{(1-\ttt)^{2d+1}}} = {\frac{h(\ttt)}{(1-\ttt)^{2d+1}}},
\]
where $h_{B_i}(\ttt)$ are the $h$--polynomials of the $B_i$. More precisely, each of the terms of $\mathcal{Z}^*$ is a $\BB$--module of the form $A\otimes \BB(-r)$ where $A$ is generated in a same degree. Such modules are isomorphic to their associated  graded modules. 

\medskip

The multiplicity of $\SS^*$ is given by the standard formula 
\[\deg \SS^* = (-1)^d   {\frac{h^{(d)}(1)}{ d!}}.\] 
We now indicate how the $h_{B_i}(\ttt)$ are  assembled. Let us illustrate
the case when $d =4$ and $i=1$. $B_1$ has  a free resolution of the strand of the Koszul complex
\[ 0 \rar \RR(-3n) \rar \RR^{5}(-2n) \lar \RR^{10}(-n) \lar \RR^{10} \lar B_1 \rar 0,\]
so that 
\[ h_{B_1}(\ttt) = 10  - 10\ttt^{n} + 5\ttt^{2n}-\ttt^{3n},\]
and similarly for all $B_i$.

\medskip

We are now ready to make our key observation. Consider a complete intersection $P$
generated by $d+1$ forms of degree $n$ in a polynomial ring of dimension $d+1$
and set $\SS^{**} = \Symi(P)$. The corresponding approximation complex now has $B_1=Z_1$. 
 The approach above would for the new $Z_i$ give the
same $h$--polynomials of the $B_i$ in the case of an almost complete intersection (but in dimension $d$).
This means that the Hilbert series of $\SS^{**}$ is given by
\[ H_{\SS^{**}}(\ttt) = {\frac{h(\ttt)}{(1-\ttt)^{2d+2}}}.\] 
 It follows 
 that
$\deg \SS^*$ can be computed as the degree of the symmetric algebra generated by a regular sequence of $d+1$ forms
of degree $n$, a  result that is given in \cite{HTU}. Thus, 
\[ \deg \SS^* = \deg \SS^{**} = \sum_{j=0}^d n^j,\]
and the calculation of $\deg \SS$ is complete. \QED

\bigskip

We will now write Theorem~\ref{degSymi} in a more convenient formulation for applications.

\begin{Theorem}\label{degSymibis} Let $\RR =k[x_1, \ldots, x_d]$ and  $I = (f_1, \ldots, f_d, f_{d+1})$ 
is an ideal  of forms of degree $n$. If $J = (f_1, \ldots, f_d)$ is a complete intersection, then
\[ \deg \SS = \sum_{j=0}^{d-1} n^j + \lambda(\RR/J:I).\]
\end{Theorem}

\demo The degree formula gives
\begin{eqnarray*}
 \deg \SS & = & \sum_{j=0}^{d-1}n^j + [n^d - \lambda(R/I)] = \sum_{j=0}^{d-1}n^j + [\lambda(\RR/J) - \lambda(\RR/I)]\\
 &=& \sum_{j=0}^{d-1} n^j +\lambda(I/J) = \sum_{j=0}^{d-1}n^j + \lambda(\RR/J:I).
 \end{eqnarray*}

\begin{Corollary} \label{degT} Let $I=(J,a)$ be an ideal of finite colength as above. Then the module of 
linear relations satisfies
 $\deg T = \length (I/J)$. In particular if $\RR[It]$ is almost Cohen--Macaulay, $T$ can be generated by $\lambda(I/J)$ elements.
\end{Corollary}

\demo From the sequence of modules of the same dimension
\[ 0 \rar T \lar \SS \lar \RR[It]\rar 0\]
we have
\[ \deg T = \deg \SS - \deg \RR[It] = \lambda(I/J).\]
\QED

The last assertion of this Corollary   can also  be obtained
from \cite[Theorem 4.1]{MPV12}.

\subsubsection*{The Cohen--Macaulay type of the module of nonlinear relations}

We recall the terminology of Cohen--Macaulay type of a module.
 Set $\BB=\RR[\TT_1, \ldots, \TT_{d+1}]$. If $E$ is a
  finitely generated $\BB$--module of codimension $r$, we say that $\Ext_{\BB}^r(E, \BB)$ is its canonical module. It is the first non vanishing $\Ext_{\BB}^i(E, \BB)$ module denoted by $\omega_E$. The minimal number of the generators of $\omega_{E}$ is called   the {\em Cohen--Macaulay type} of $E$ and is denoted by  $\mbox{\rm type}(E)$.
When $E$ is graded and Cohen--Macaulay, it gives the last Betti number of a projective resolution of $E$. It can be expressed in different ways, for example for
the module of linear relations $\omega_T=\Ext_{\BB}^d(T, \BB) = \Hom_{\SS}(T,\omega_{\SS})$.

\begin{Proposition}\label{typeofT} Let $\RR$ be a Gorenstein local  ring of dimension $d \geq 2$
 and $I=(J,a)$ an ideal of finite colength as above. If $\RR[It]$ is an aCM algebra and $\omega_{R[It]}$ is Cohen--Macaulay, then the type of the  module $T$ of nonlinear relations satisfies
\[ \mbox{\rm type}(T) \leq \mbox{\rm type }(S_J(I)) + d-1,\]
where $S_J(I)$ is the Sally module.
\end{Proposition}

\demo We set $\Rees = \RR[It]$ and $\Rees_0 = \RR[Jt]$.
First apply $\Hom_{\BB}(\cdot, \BB) $ to the basic presentation
\[ 0 \rar T \lar \SS \lar \Rees \rar 0,\]
to obtain the cohomology sequence
\begin{eqnarray}\label{type1}
 0 \rar \omega_{\Rees} \lar \omega_{\SS} \lar \omega_T \lar \Ext_{\BB}^{d+1}(\Rees, \BB) \rar 0.
 \end{eqnarray}

Now apply the same functor to the exact sequence of $\BB$--modules
\[ 0 \rar I\cdot \Rees[-1]  \lar \Rees \lar \RR \rar 0\] 
to obtain the exact sequence
\[ 0 \rar \omega_{\Rees} \stackrel{\theta}{\lar} \omega_{I\Rees[-1]} \lar \Ext_{\BB}^{d+1}(\RR, \BB) = \RR
  \lar \Ext_{\BB}^{d+1}(\Rees, \BB) \lar
\Ext_{\BB}^{d+1}(I\Rees[-1], \BB) \rar 0.\]  
Since $\omega_{\Rees}$ is Cohen--Macaulay and $\dim \RR \geq 2$,
 the cokernel of $\theta$ is either $\RR$ or an $\m$-primary ideal that satisfies the condition $S_2$ of Serre. The only choice is $\coker(\theta)=\RR$. Therefore 
 \[\Ext_{\BB}^{d+1}(\Rees, \BB)\simeq \Ext_{\BB}^{d+1}(I\Rees[-1], \BB).\]
 
 Now we approach the module $ \Ext_{\BB}^{d+1}(I\Rees, \BB)$ from a different direction. We note that $\Rees$---but not $\SS$ and $T$---is also a finitely generated $\BB_0=\RR[\TT_1, \ldots, \TT_d]$--module as it is annihilated by a monic polynomial $\ff$ in $\TT_{d+1}$ with coefficients in $\BB_0$. By Rees Theorem
  we have  that for all $i$, $\Ext_{\BB}^i(\Rees, \BB) =\Ext_{\BB/(\ff)}^{i-1}(\Rees, \BB/(\ff))$, and a similar observation applies to $I\cdot \Rees$.
  
  \medskip
  
  Next consider the finite, flat morphism $\BB_0 \rar \BB/(\ff)$. For any $\BB/(\ff)$--module $E$ with a projective
  resolution $\mathbb{P}$, we have that $\mathbb{P}$ is a projective $\BB_0$--resolution of $E$. This 
  means that the isomorphism of complexes
  \[ \Hom_{\BB_0}(\mathbb{P}, \BB_0) \simeq \Hom_{\BB/(\ff)}(\mathbb{P}, \Hom_{\BB_0}(\BB/(\ff), \BB_0))=
 \Hom_{\BB/(\ff)}(\mathbb{P}, \BB/(\ff)) 
  \]
  gives  isomorphisms for all $i$
  \[ \Ext_{\BB_0}^i(E, \BB_0) \simeq  \Ext_{\BB/(\ff)}^i(E, \BB/(\ff)). \]  
 Thus 
 \[ \Ext_{\BB}^i(\Rees, \BB) \simeq \Ext_{\BB_0}^{i-1}(\Rees, \BB_0).\]  
In particular, $\omega_{\Rees} = \Ext_{\BB_0}^{d-1}(\Rees, \BB_0)$.
 

\medskip

Finally apply $\Hom_{\BB_0}(\cdot, \BB_0)$ to the exact sequence of $\BB_0$--modules
 and examine its   cohomology sequence.

\[ 0 \rar I\cdot \Rees_0 \lar I\cdot \Rees \lar S_J(I) \rar 0\]
is then
\[ 0 \rar \omega_{I\Rees} \lar \omega_{I\Rees_0} \lar \omega_{S_J(I)}\lar \Ext_{\BB_0}^d(\Rees, \BB_0)= \Ext_{\BB}^{d+1}(\Rees, \BB) \rar 0.\] 
 Taking this into
(\ref{type1}) and the  that $\mbox{\rm type}(\SS)=d-1$ gives the desired estimate.
 \QED

\medskip

\begin{Remark}{\rm
A class of ideals with $\omega_{\Rees}$ Cohen--Macaulay is discussed in Corollary~\ref{canofRees}(b).
}\end{Remark}

\section{Distinguished aCM algebras}

This section treats several classes of Rees algebras which are almost Cohen--Macaulay.

\subsection{Equi-homogeneous acis}

We shall now treat an important class of extremal Rees algebras.
Let $\RR = k[x_1, \ldots, x_d]$ and let $I=(a_1, \ldots, a_d, a_{d+1})$ be an ideal
of finite colength, that is,  $\m$--primary. We further assume that the first
 $d$ generators form a regular sequence and $a_{d+1} \notin  (a_1, \ldots, a_d)$. 
    If $\deg a_i=n$, the integral closure
 of $J=(a_1, \ldots, a_d)$ is the ideal $\m^n$, in particular $J$ is a minimal reduction of $I$. The integer
 $\edeg(I) = \red_J(I) +1$ is called the {\em elimination degree} of $I$.
 The study of the equations of $I$, that is, of $\RR[It]$, depends on a comparison between
 the metrics  of $\RR[It]$ to
 those of $\RR[\m^n t]$, which are well known.

\begin{Proposition} {\rm (\cite{syl2})} \label{birideal}
The following conditions are equivalent:
\begin{itemize}
\item[{\rm (i)}] $\Phi$ is a birational mapping, that is the natural embedding $\mathcal{F}(I) \hookrightarrow \mathcal{F}(\m^n)$ is an isomorphism of quotient fields;

\item[{\rm (ii)}] $\red_J(I) = n^{d-1}-1$;

\item[{\rm (iii)}] $e_1(I) = {\frac{d-1}{2}}(n^d - n^{d-1})$;

\item[{\rm (iv)}] $\RR[It]$ satisfies the condition $R_1$ of Serre.
\end{itemize}
\end{Proposition}

For lack of a standardized terminology, we say that these ideals are {\em birational}\label{birational ideal}.

\begin{Corollary}\label{canofRees} For an ideal $I$ as above, the following hold:
\begin{itemize}

\item[{\rm (i)}] The algebra $\RR[It]$ is not Cohen--Macaulay except when $I = (x_1, x_2)^2$.

\item[{\rm (ii)}] The canonical module of $\RR[It]$ is Cohen--Macaulay.

\end{itemize}

\end{Corollary}

\demo (i) follows from the condition of Goto--Shimoda (\cite{GS82}) that the reduction number of a Cohen--Macaulay Rees algebra $\RR[It]$ must satisfy $\red_J(I) \leq \dim \RR -1$, which in the case 
$n^{d-1}- 1 \leq d-1$
is only met if $d =n =2$ 

\medskip

\noindent (ii) The embedding $\RR[It] \hookrightarrow \RR[\m^n t]$ being an isomorphism in codimension one, their canonical modules are isomorphic. The canonical module of a Veronese subring such as
$\RR[\m^n t]$ is well-known (see \cite[p. 187]{HV85}, \cite{HSV87}; see also \cite[Proposition 2.2]{BR}).

\subsection*{Binary ideals}

These are the ideals of $\RR=k[x,y]$ generated by $3$ forms of degree $n$. Many of their Rees algebras are almost Cohen--Macaulay. We will showcase the technology of the Sally module in treating a much studied class of ideals. First
we discuss a simple case (see also \cite{syl1})).

\begin{Proposition}\label{22} Let $\phi$ be a $3\times 2$ matrix of quadratic forms in $\RR$ and $I$ the ideal
given by its $2\times 2$ minors. Then $\RR[It]$ is almost Cohen--Macaulay.
\end{Proposition}

\demo These ideals have reduction number $1$ or $3$. In the first case all of its Sally modules vanish and 
$\RR[It]$ is Cohen--Macaulay.

\medskip

In the other case $I$ is a birational ideal and $e_1(I) = {4\choose 2} = 6$. A simple calculation shows that
$\lambda(\RR/I) = 12$, so that $\lambda(I/J) = 16-12 = 4$. To apply Theorem~\ref{Huckaba}, we need to verify
the equation 
\begin{eqnarray}\label{Sally22}
 f_1 + f_2 + f_3 = 6.
 \end{eqnarray}
We already have $f_1=4$. To calculate $f_2$ we need to take $\lambda(R/I_1(\phi))$ in Corollary~\ref{F2}. $I_1(\phi)$ is an ideal generated by $2$ generators or $I_1(\phi) = (x,y)^2$. But in the first case the Sylvester resultant
of
 the linear equations of $\RR[It]$ would be a quadratic polynomial, that is $I$ would have reduction number $1$, which 
 would contradict the assumption. Thus by Corollary~\ref{F2}, $f_2 = 4-\lambda(\RR/I_1(\phi)) = 1$. Since
 $f_2\geq f_3>0$ we have $f_3=1$ and the equation (\ref{Sally22}) is satisfied. \QED

\bigskip

We have examined higher degrees examples of birational ideals of this type which are/are not almost Cohen--Macaulay.
Quite a lot is known about the following ideals.
$\RR = k[x,y]$ and $I$ is a codimension $2$ ideal given by that $2\times 2$ minors of a $3\times 2$ matrix
with homogeneous entries of degrees $1$ and $n-1$. 

\begin{Theorem} \label{2birideal}
If $I_1(\phi) = (x, y)$ then:
\begin{enumerate} 

\item[{\rm (i)}] $\deg \mathcal{F}(I) = n$, that is $I$ is birational.

\medskip

\item[{\rm (ii)}]  $\RR[It]$ is almost Cohen--Macaulay.

\medskip

\item[{\rm (iii)}]  The equations of $\LL$ are given by a straightforward algorithm.

\end{enumerate}
\end{Theorem}

\demo The proof of (i) is in \cite{CHW},  and in other sources (\cite{CdA}, \cite{KPU}; see also
\cite[Theorem 2.2]{DHS} for a broader statement in any characteristic and \cite[Theorem 4.1]{Simis04} in
characteristic zero),
 and of  (ii)  in \cite[Theorem 4.4]{KPU}, while (iii)  was conjecturally  given in  \cite[Conjecture 4.8]{syl1} and proved in \cite{CHW}. We give a combinatorial proof of (ii) below (Proposition~\ref{aCMofbin}).
\QED

\medskip

We note that
$\deg(\SS) = 2n$, since $\SS$ is a complete intersection defined by two forms of [total] degrees $2$ and $n$, while $\RR[Jt]$ is defined by one equation of degree $n+1$. Thus $\nu(T) \leq 2n-(n+1)= n-1$, which is the number of generators given in the algorithm.

\bigskip

We point out a property of the module $T$. We recall that an $\AA$--module is an {\em Ulrich} module if it is a maximal Cohen--Macaulay module with $\deg M=\nu(M)$ (\cite{HK}). 

\begin{Corollary} $T$ is an Ulrich $\SS$--module.
\end{Corollary}


Considerable numerical information in the Theorem~\ref{2birideal} will follow from:

\begin{Proposition}\label{aCMofbin} If $\deg \alpha =1$ and $\deg \beta = n-1$, then $\length(F_j) = n -j$. In particular, $\RR[It]$ is almost Cohen--Macaulay.
\end{Proposition}

\demo
Note that  the ideal is birational, $F_{n-1}\neq 0$. On the other hand, $\LL$ contains fresh generators in all degrees $j\leq n$. This means that for $f_j = \lambda(F_j)$, 
\[ f_j> f_{j+1}>0, \quad j<n.\]
Since $f_1 = n-1$, the decreasing sequence of integers 
\[ n-1 = f_1 > f_2 > \cdots > f_{n-2} > f_{n-1}> 0\]
implies that $f_j = n-j$.
Finally, applying Theorem~\ref{Huckaba} we have that $\RR[It]$ is an aCM algebra since $\sum_{j}f_j = e_1(I)= {n\choose 2}$. \QED

\subsection*{Quadrics}

Here we explore sporadic classes of aCM algebras defined by   quadrics in
$k[x_1, x_2, x_3, x_4]$.

\medskip

  First we
 use Proposition~\ref{F2} to look at  other cases of quadrics. For $d=3$, $n=2$, $\edeg(I)=2$ or $4$. In the first case
$J:a= I_1(\phi)$. In addition $J:a\neq \m$ since the socle degree of
$\RR/J$ is $3$. Then $\RR[It]$ is Cohen--Macaulay. If $\edeg(I)=4$ we must have [and conversely!]
 $\length(\RR/J:a)=2$ and $I_1(\phi)=\m$. Then
$\RR[It]$ is almost Cohen--Macaulay.

\medskip

Next we treat almost complete intersections of finite colength generated by quadrics of $\RR= k[x_1, x_2, x_3, x_4]$. Sometimes we denote the variables by $x,y,\ldots$, or use these symbols to denote  [independent] linear forms. For notation we
use $J= (a_1, a_2, a_3, a_4)$, and $I = (J, a)$.

\bigskip

Our goal  is to address the following:
\begin{Question}{\rm  Let $I$ be an almost complete intersection generated by $5$ quadrics of $x_1, x_2, x_3, x_4$. If   
  $I$ is a birational ideal, in which cases   is $\RR[It]$ is an almost Cohen--Macaulay algebra? In this case, what are the generators of the its module of nonlinear relations?
}\end{Question} 

In order to make use of Theorem~\ref{Huckaba}, our main tools are Corollary~\ref{F2} and \cite[Theorem 2.2]{syl2}. 
They make extensive use of the syzygies of $I$. The question forks into three cases, but our analysis is 
complete in only one of them.

\subsubsection*{The Hilbert functions of quaternary quadrics}

We make a quick classification of the Hilbert functions of the ideals $I=(J,a)$. Since $I/J \simeq \RR/J:a$ and $J$ is a complete intersection, the problem is equivalent to determine the Hilbert functions of $\RR/J:a$, with $J:a$ a Gorenstein ideal.
The Hilbert function $H(\RR/I)$ of $\RR/I$ is $H(\RR/J)-H(\RR/J:a)$. We will need the Hilbert function of the corresponding canonical module in order to make use of \cite[Proposition 3.7]{syl2} giving information about $L_2/\BB_1L_1$.

\medskip

We shall refers to the  sequence   $(f_1, f_2, f_3, \ldots, )$, $f_i = \lambda(I^i/JI^{i-1})$, as the $\ff$--sequence of     $(I,J)$. We recall that these sequences are monotonic and that 
 if $I$ is birational, $\sum_{i\geq  1} f_i=e_1(I)=12$. 

\begin{Proposition}
Let  $\RR= k[x_1, x_2, x_3, x_4]$ and $I=(J,a)$ an almost complete intersection generated by $5$ quadrics, where $J$ is a complete intersection,.  Then 
$L=\lambda(\RR/J:a)\leq 6$ and the possible Hilbert functions of $\RR/J:a)$ are:
\begin{eqnarray*}
L = 6 & : & (1,4,1), \quad (1,2,2,1)^{*}, \quad (1,1,1,1,1,1)^{*} \\
L = 5 & : & (1,3,1), \quad (1,1,1,1,1)^{*} \\
L = 4 & : & (1,2,1), \	\quad (1,1,1,1)^{*} \\
L = 3 & : & (1,1,1)^{*} \\
L = 2 & : & (1,1)^{**} \\
L = 1 & : & (1)^{**}
\end{eqnarray*}
If $I$ is a birational ideal, the corresponding Hilbert function is one of the unmarked sequences above.
\end{Proposition}

\demo Since $\lambda(\m^2/I)\geq 5$ and $\lambda(\m^2/J)=11$, $L=\lambda(\RR/J:a)\leq 6$. Because the Hilbert function of $R/(J:a)$ is symmetric and $L \leq 6$, the list includes all the viable Hilbert functions.

\medskip

Let us first rule out those marked with ${\mbox{\rm a}}^{*}$, while those marked with ${\mbox{\rm a}}^{**}$ cannot be birational.
In each of these $J:a$ contains at least $2$ linearly independent linear forms, which we denote by $x,y$, so that $J:a/(x,y)$ is a Gorenstein ideal of the regular ring $\RR/(x,y)$. It follow that $(J:a)/(x,y)$ is a complete intersection. In the case of $(1,2,2,1)$, $J:a = (x,y, \alpha, \beta)$, where $\alpha$ is a form of degree $2$ and $\beta$ a form of degree $3$, since $\lambda(\RR/J:a)=6$. Since $J\subset J:a$, all the generators of $J$ must be contained in $(x,y,\alpha)$, which is impossible by Krull theorem.
Those  strings with at least  $3$ $1$'s are also excluded since $J:a$ would
 have the form $(x,y,z,w^s)$, $s\geq 3$, and the argument above applies. The case $(1,1)$, $J:a = (x,y,z,w^2)$. This means that $I_1(\phi)=J:a$, or $J:a = \m$. In the first case, by Corollary~\ref{F2Cor}, $I^2 = JI$. In the second case, $I_1(\phi) = \m$. This will imply that
$ \lambda(I^2/JI) = 2 - 1 = 1$, and therefore $I$ will not be birational (need the summation to add to $12$).
\QED

\subsubsection*{Hilbert function $(1,4,1)$}

If $R/J:a$ has Hilbert function $(1,4,1)$, $J:a \subset \m^2$ but we cannot have equality since $\m^2$ is not a Gorenstein ideal. We also have
 $I_1(\phi) \subset \m^2$. If they are not equal, $I_1(\phi) = J:a$, which by Corollary~\ref{F2Cor} would mean that $\red_J(I) = 1$.

\begin{Theorem} \label{141} Suppose $I$  that is birational and $I_1(\phi)\subset \m^2$. 
Then
$\RR[It]$ is almost Cohen--Macaulay.
\end{Theorem}

\demo
The assumption $I_1(\phi)\subset \m^2$ means that the Hilbert function of $J:a$ is $(1,4,1)$ and vice-versa.
 Note also that by assumption  $\length(I^7/JI^6)\neq 0$. Since $\lambda(I/J) = \lambda(\RR/J:a) = 6$, it suffices to show that $\lambda(I^2/JI)=1$.
 From $\lambda(\m^2/I) = 5$, the module $\m^2/I$ is of length $5$ minimally generated by $5$ elements. Therefore $\m^3\subset I$, actually $\m^3=\m I$.

\medskip
 
There is an isomorphism  $\hh: \RR/J:a \simeq I/J$, $r\mapsto ra$. It moves the socle of 
$\RR/J:a$ into the socle of $I/J$. If $a\notin J:a$, then $\m^2 = (J:a, a)$ and $a$ gives the
socle of $\RR/J:a$, thus it is mapped by $\hh$ into the socle of
$I/J$, that is $\m\cdot a^2\in J$. Thus,
$ \m\cdot a^2\in \m^3 J \subset JI$. On the other hand, if  $a\in J:a$, then, since $a^2\in J$, we have $a^2\in \m^2 J$ and 
$\m\cdot a^2\in \m^3 J\subset JI$.

\medskip

An example is $J= (x^2, y^2, z^2, w^2)$, $a =  xy + xz + xw + yz$. 
 \QED

\subsubsection*{Hilbert function $(1,3,1)$}

Our discussion about this case is very spare.

\begin{itemize}

\item
For these Hilbert functions, 
$ J:a = (x, P)$, 
where $P$ is a Gorenstein ideal in a regular local ring of dimension $3$--and therefore 
is given by the Pfaffians of a 
skew-symmetric matrix, necessarily $5\times 5$.
 Since $J \subset J:a$, $L$ must contain forms of degree $2$. In addition,   $P$ is given by $5$  $2$--forms (and $(x,\m^2)/(x,P)$ is the socle of $\RR/J:a$).
  
\item If $I$ is birational,  $\RR[It]$ is almost Cohen--Macaulay if and only if $\lambda(I^2/JI) = 2$ and $\lambda(I^3/JI^2)=1$.   The first equality, by Proposition~\ref{F2}, requires $\lambda(\RR/I_1(\phi)) = 3$ which gives that $I_1(\phi)$ contains the socle of $J:a$ and another independent linear form. In all it  means that $I_1(\phi) = (x,y, (z,w)^2)$. On the other hand
 $\lambda(I^2/JI)=2$ means that $JI:I^2 = (x,y,z,w^2)$ (after more label changes).

\item  An   example is  $J=(x^2, y^2, z^2, w^2)$ with $a=  xy + yz + zw + wx + yw$. The ideal $I= (J,a)$ is birational.

 \end{itemize}
  \subsubsection*{Hilbert function $(1,2,1)$}
 
 We do not have the full analysis of this case either.
 
 \begin{itemize}
 
 \item  An example is  $J=(x^2, y^2, z^2, w^2)$ and $a = xy+yz+xw+zw$. The ideal $I=(J,a)$  is birational.
The expected $\ff$--sequence of such  ideals is $(4,3,1,1,1,1,1)$.

\item
If $I$ is birational then $I_1(\phi)=\m$. We know that $I_1(\phi)\neq J:a$, so $I_1(\phi) \supset \m^2$, that is
$I_1(\phi)= (x,y,\m^2)$,  $(x,y,z,\m^2)$, or $\m$. Let us exclude the first two cases.

$(x,y,\m^2)$: This leads to two equations
\begin{eqnarray*}
xa &=& xb + yc\\
ya &=& xd + ye,
\end{eqnarray*} 
with $b,c,d,e\in J$. But this gives the equation $(a-b)(a-e) -dc=0$, and $\red_J(I)\leq 1$.

\medskip

  $(x,y,z,\m^2)$: Then the Hilbert function of $\RR/I_1(\phi)$ is $(1,1)$. According to \cite[Proposition 3.7]{syl2}, $L_2$ has a form of bidegree $(1,2)$, with coefficients in $I_1(\phi)$, that is, in
  $(x,y,z)$. This  gives $3$ forms with coefficients in this ideal, two in degree $1$, so by elimination we get a monic
 equation of degree $4$.

\end{itemize}

We summarize the main points of these observations into a normal form assertion.

\begin{Proposition}\label{121} Let $I$ be a birational ideal and the  Hilbert function of $\RR/J:a$ is
$(1,2,1)$. Then up to a change of variables 
to $\{ x,y,z,w\}$, $I$ is a Northcott ideal, that is there is a $4\times 4$ matrix $\AA$,
\[ \AA = \left[ \begin{array}{ccc}
& \BB &  \\
 \hline 
 & \mathbf{C} & \\
 \end{array} \right] . \]
 where $\BB$ is a $2\times 4$ matrix whose entries are linear forms and $\mathbf{C}$ is a matrix with scalar entries and 
 $\mathbf{V} = [x,y, \alpha, \beta]$, where $\alpha, \beta$ are quadratic forms in $z,w$ such that
 \[ I = (\mathbf{V}\cdot \AA, \det \AA).\] 
 \end{Proposition}
 
 \demo There are two independent linear forms in $J:a$ which we denote by $x,y$.   
 We observe that $(J:a)/(x,y)$ is a Gorenstein ideal in a polynomial ring of  
  dimension two, so it is a complete intersection: $J:a = (x,y, \alpha,\beta)$, with $\alpha$ and $\beta$ forms of degree $2$ (as $\lambda(\RR/J:a)=4$), from which we remove the terms in $x,y$, that is we may assume $\alpha, \beta\in (z,w)^2$.
  
  \medskip
  
  Since $J \subset J:a$, we have a matrix $\AA$,
 \[ J = [x,y,\alpha, \beta] \cdot \AA =\mathbf{V} \cdot \AA.\]
 By duality $I=J:(J:a)$, which by Northcott theorem (\cite{Northcott}) gives 
 \[ I= (J, \det \AA).\]
 Note that $a$ gets, possibly, replaced by $\det \AA$.
  The statement about the degrees of the entries of $\AA$ is clear. \QED

 \begin{Example}{\rm
 Let 
 \[ \AA = \left[ \begin{array}{rrrr}
 x+y & z + w & x-w & z \\
 z & y+w & x - z & y \\
 1 & 0 & 2 & 3 \\
 0 & 1 & 1 & 2 \\
 \end{array} \right], \quad {\mathbf v} = \left[ x, y,  z^2 + zw + w^2, z^2-w^2\right] .\]
This  ideal is birational but $\RR[It]$ is not aCM. This is unfortunate but opens the question of when such ideals are birational. The $\ff$--sequence here is $(4,3,3,1,1,1,1)$.

 }\end{Example}

\subsubsection*{The degrees of $\LL$} 

We examine how the Hilbert function of $\RR/J:a$ organizes
 the generators of $\LL$. We denote the presentations variables by
$u, \TT_1, \TT_2, \TT_3, \TT_4$, with $u$ corresponding to $a$. 

\begin{itemize}

\item $(1,4,1)$: We know (Theorem~\ref{141}) that $JI:I^2= \m$. This means that we have forms 
\begin{eqnarray*}
\hh_1  & = & xu^2 + \cdots \\
\hh_2  & = & yu^2 + \cdots \\
\hh_3 & = & zu^2 + \cdots \\
\hh_4 & = & w u^2 + \cdots 
\end{eqnarray*}
with the $(\cdots)$ in $(\TT_1, \TT_2, \TT_3, \TT_4)\BB_1$. The corresponding resultant, of degree $8$, is nonzero. 

\item $(1,2,1)$: There are two forms of degree $1$ in $\LL$,
\begin{eqnarray*}
\ff_1 & = & x u + \cdots\\
\ff_2 & = & y u + \cdots 
\end{eqnarray*}
The forms in $L_2/\BB_1L_1$ have coefficients in $\m^2$. This will follow from $I_1(\phi)= \m$. We need a way to generate
two forms of degree $3$. Since we expect $JI:I^2 = \m$, this would mean the presence of two forms in $L_3$,
\begin{eqnarray*}
\hh_1^{*} & = & z u^3 + \cdots\\
\hh_2^{*} & = & w u^3 + \cdots ,
\end{eqnarray*}
which together with $\ff_1$ and $\ff_2$ would give the nonzero degree $8$ resultant.%

\item $(1,3,1)$: There is a form $\ff_1 = xu + \cdots \in L_1$ and two forms in $L_2$
\begin{eqnarray*}
\hh_1 & = & y u^2 + \cdots \\
\hh_2 & = & zu^2+ \cdots 
\end{eqnarray*}
predicted by \cite[Prop 3.7]{syl2}  if $I_1(\phi) = (x,y,z,w^2)$. (There are indications that this is always the case.) 
 We need a cubic equation $\hh_3^{*} = wu^3+ \cdots$ to launch the nonzero resultant of degree $8$.

\end{itemize}

For all quaternary quadrics with $\RR[It]$ almost Cohen--Macaulay, Corollary~\ref{degT} says that 
$\nu(T) \leq \lambda(\RR/J:a)$. Let us compare to the actual number of generators in the examples discussed above:

\[  \left[
\begin{array}{ccc}
\nu(T) &  & \lambda(I/J)\\
5 & (1,4,1) & 6\\
4 & (1,3,1) & 5\\
4 & (1,2,1)& 4
\end{array} \right] 
\]
We note that in the last case, $T$ is an  Ulrich module.

\subsection{Monomial ideals}

Monomial ideals of finite colength which are almost complete intersections have a very simple description. We examine a narrow class of them.
Let $\RR=k[x,y,z]$ be  a polynomial ring over and infinite field and let $J$ and $I$ be $\RR$--ideals such that
\[{\ds J=(x^{a},\; y^{b},\; z^{c}) \subset (J,\; x^{\alpha} y^{\beta} z^{\gamma})=I. }\] 
This is the general form of almost complete intersections of $\RR$ generated by monomials. Perhaps the most interesting cases are those 
where
${\ds \frac{\alpha}{a} + \frac{\beta}{b} + \frac{\gamma}{c} <1}$. This inequality
ensures that $J$ is not a reduction of $I$.
Let 
\[{\ds Q =(x^{a}-z^{c},\; y^{b}-z^{c},\;x^{\alpha} y^{\beta} z^{\gamma} ) }\]
 and suppose that ${\ds a > 3 \alpha,\; b > 3 \beta,\;
c > 3 \gamma }$. Note that $I=(Q,\; z^c)$.  Then $Q$ is a minimal reduction of $I$ and the reduction number $\mbox{\rm red}_{Q}(I) \leq
2$. We label these ideals $I(a,b,c,\alpha, \beta, \gamma)$.

\medskip

We will examine in detail the case $a=b=c=n\geq 3$ and
$\alpha=\beta = \gamma =1$. We want to argue that
$\RR[It]$ is almost Cohen--Macaulay. To benefit from the monomial generators in using {\em Macaulay2} we set
$I = (xyz, x^n, y^n,z^n)$. Setting $\BB=\RR[u, \TT_1, \TT_2, \TT_3]$, we claim that
\[
\LL= (z^{n-1}u - xy\TT_3, y^{n-1}u - xz\TT_2, x^{n-1}u - yz\TT_1, z^n\TT_2-y^n\TT_3,  z^n\TT_1-x^n\TT_3, 
y^n\TT_1-x^n\TT_2, 
\]
\[ y^{n-2}z^{n-2}u^2 - x^2\TT_2\TT_3, x^{n-2} z^{n-2} u^2 -y^2 \TT_1\TT_3,   x^{n-2}y^{n-2} u^2 - z^2\TT_1\TT_2,
x^{n-3}y^{n-3}z^{n-3}u^3 - \TT_1\TT_2\TT_3). 
  \]
 We also want to show that these ideals define an almost Cohen--Macaulay Rees algebra.

\bigskip

There is a natural specialization argument. Let $X$, $Y$ and $Z$ be new indeterminates and let
$\BB_0 = \BB[X,Y,Z]$. In this ring define the ideal $\LL_0$ obtained by replacing in the list above of generators 
of $\LL$, $x^{n-3}$ by $X$ and accordingly $x^{n-2}$ by $xX$, and so on; carry  out similar actions on the other variables:

\[
\LL_0= (z^2 Zu - xy\TT_3, y^{2}Yu - xz\TT_2, x^{2}Xu - yz\TT_1, z^3Z\TT_2-y^3Y\TT_3,  z^3Z\TT_1-x^3X\TT_3, 
y^3 Y\TT_1-x^3X\TT_2, 
\]
\[ yzYZu^2 - x^2\TT_2\TT_3, x z XZ u^2 -y^2 \TT_1\TT_3,   xyXY u^2 - z^2\TT_1\TT_2,
XYZu^3 - \TT_1\TT_2\TT_3). 
  \]

Invoking {\em Macaulay2} gives a (non-minimal) projective resolution
\[ 0 \rar \BB_0^4 \stackrel{\phi_4}{\lar}
 \BB_0^{17} \stackrel{\phi_3}{\lar}
\BB_0^{22} \stackrel{\phi_2}{\lar}
 \BB_0^{10} \stackrel{\phi_1}{\lar}
  \BB_0 \lar \BB_0/\LL_0 \rar 0.
 \]

We claim that the specialization $X \rar x^{n-3}$, $Y \rar y^{n-3}$, $Z \rar z^{n-3}$ gives a projective resolution
of $\LL$.

\begin{itemize}

\item Call $\LL'$ the result of the specialization in $\BB$. We argue that $\LL' = \LL$.

\medskip

\item Inspection of the Fitting ideal $F$ of  $\phi_4$ shows that it contains $(x^3, y^3,z^3, u^3, \TT_1\TT_2\TT_3)$. 
From standard theory, the radicals of the Fitting ideals of $\phi_2$ and $\phi_2$  contain $\LL_0$, and therefore
the radicals of the Fitting ideals of these mappings after specialization will contain the ideal $(L_1)$ of $\BB$, as 
$L_1 \subset \LL'$.

\medskip

\item Because $(L_1)$ has codimension $3$,
 by the acyclicity 
 theorem (\cite[1.4.13]{BH}) 
  the complex gives a projective resolution of $\LL'$. Furthermore, as $\mbox{\rm proj. dim }\BB/\LL' \leq 4$,
 $\LL'$ has no associated primes of codimension $\geq 5$. Meanwhile
the Fitting ideal of $\phi_4$ having codimension $\geq 5$,  forbids the existence of associated primes
of codimension $4$. Thus $\LL'$ is 
unmixed.

\medskip

\item Finally, in $(L_1) \subset \LL' $, as $\LL'$ is unmixed its associated primes are minimal primes of 
$(L_1)$, but by Proposition~\ref{canoseq}(iii), there are just two such, $\m\BB$ and $\LL$. Since
$\LL' \not \subset \m\BB$, $\LL$ is its unique associated prime. Localizing at $\LL$ gives
the equality of $\LL' $ and $\LL$ since $\LL$ is a primary component of $(L_1)$.

\end{itemize}

Let us sum up this discussion:

\begin{Proposition} \label{nnn111} The Rees algebra of $I(n, n, n, 1, 1, 1)$, $n\geq 3$, is almost Cohen--Macaulay.
\end{Proposition}

\begin{Corollary} $e_1(I(n,n,n,1,1,1)) = 3(n+1)$.
\end{Corollary}

\demo Follows easily since $e_0(I) = 3n^2$, the colengths of the monomial ideals $I$ and $I_1(\phi)$ directly
calculated and $\red_J(I) = 2$ so that 
\[ e_1(I) =  \lambda(I/J) + \lambda(I^2/JI) = \lambda(I/J) + [\lambda(I/J) - \lambda(\RR/I_1(\phi))]=
(3n-1) + 4.\]
\QED

\begin{Remark}{\rm 
We have also experimented with other cases beyond those with $xyz$ and in higher dimension as well. 
\begin{itemize}

\item In $\dim \RR=4$, the ideal $I= I(n,n,n,n,1,1,1,1)= (x_1^n, x_2^n, x_3^n, x_4^n, x_1x_2x_3x_4)$, $n\geq 4$, has
a Rees algebra $\RR[It]$ which is almost Cohen--Macaulay.

\item The argument used was a copy of the previous case, but we needed to make an adjustment in the last step
to estimate the codimension of the last Fitting ideal $F$ of the corresponding mapping $\phi_5$. This is a large matrix,
so it would not be possible to find the codimension of $F$ by looking at all its maximal minors. Instead, one argues as follows.
Because $I$ is $\m$--primary, $\m\subset \sqrt{F}$, so we can drop the entries in $\phi_5$ in $\m$ Inspection will give
$u^{16} \in F$, so dropping all $u$'s gives an additional minors in $\TT_1, \ldots, \TT_4$, for $\height (F)\geq 6$.
This suffices to show that $\LL = \LL'$.


\end{itemize}

}\end{Remark}

\begin{Conjecture} \label{monoaciRacm}
{\rm Let $I$ be a monomial ideal of $k[x_1, \ldots, x_n]$. If $I$ is an almost complete intersection 
of finite colength its Rees algebra $\RR[It]$ is almost Cohen-Macaulay.

}\end{Conjecture}


\begin{thebibliography}{99}
\addcontentsline{toc}{section}{Bibliography}






\bibitem{BH}{W. Bruns and J. Herzog, {\em Cohen--Macaulay Rings}, Cambridge University Press, 1993.}
 

\bibitem{BR}{W. Bruns and G. Restuccia, Canonical modules of Rees algebras, {\em J. Pure \& Applied Algebra} 
{\bf 201} (2008), 189--203.}

\bibitem{CHV}{A. Corso, C. Huneke and W. V. Vasconcelos, On the integral closure of ideals, {\em Manuscripta Math.}
{\bf 95} (1998), 331--347.}


\bibitem{CPV1}{A. Corso, C. Polini and  W. V. Vasconcelos, { Links of
prime ideals}, {\em Math.  Proc. Camb.  Phil. Soc.} {\bf 115} (1994), 431--436.}


\bibitem{CdA}{T.  Cortadellas and C. D'Andrea,  Rational plane curves parametrizable by conics, arXiv:1104.2782 [math.AG].
}

\bibitem{CHW}{D. Cox, J. W. Hoffman and  H. Wang, Syzygies and the Rees
algebra, {\em J. Pure \& Applied Algebra} {212} (2008)
1787--1796.}

\bibitem{DHS}{A. V. Doria, H. Hassanzadeh and A. Simis, A characteristic free criterion of birationality,
{\em Advances in Math.} {\bf 230} (2012), 390--413.}

\bibitem{Eisenbudbook} {D. Eisenbud, {\em Commutative Algebra with a view toward Algebraic Geometry}, Springer, 1995.  }

\bibitem{ehrhart}{C. Escobar, J.
Mart\'\i nez-Bernal and R. H. Villarreal, Relative volumes and minors
 in monomial subrings, {\em Linear Algebra Appl.} {374} (2003) 275--290.}

\bibitem{GS82} {S. Goto and Y. Shimoda,  Rees algebras of Cohen--Macaulay
     local rings, in {\em Commutative Algebra}, Lect. Notes in Pure and Applied Math. {\bf 68}, Marcel Dekker,
     New York, 1982, 201--231.} 
 
\bibitem{Macaulay2}{D. Grayson and M. Stillman,
  {Macaulay 2}, a software system for research in algebraic
  geometry, 2006.
  Available at {\tt http://www.math.uiuc.edu/Macaulay2/}.}

\bibitem{HK}{J. Herzog and M. K\"{u}hl, Maximal Cohen--Macaulay modules over Gorenstein rings and Bourbaki sequences, {\em  Commutative Algebra and Combinatorics}, Adv. Stud. Pure Math. {\bf 11} (1987), 65--92.}

\bibitem{HSV1}{J. Herzog, A. Simis and W. V. Vasconcelos, Approximation complexes of blowing-up rings, 
{\em J. Algebra} {\bf 74} (1982), 466--493.}


\bibitem{HSV83}{J. Herzog, A. Simis and W. V. Vasconcelos,
Koszul homology and blowing-up rings, in Proc. 
Trento Conf. in Comm. Algebra, Lect. Notes in Pure and Applied Math. {\bf 84}, Marcel Dekker,
New York, 1983, 79--169.}

\bibitem{HSV87}{J. Herzog, A. Simis and W. V. Vasconcelos,
On the canonical module of the  Rees algebras and the associated graded ring of an ideal, 
{\em J. Algebra} {\bf 105} (1987), 205--302.} 

\bibitem{HTU}{J. Herzog, N. V. Trung and B. Ulrich,      On the multiplicity of blow-up rings of ideals 
generated by d--sequences, {\em J. Pure \& Applied Algebra} {\bf 80} (1992), 273--297.   }


\bibitem{HV85}{J. Herzog and W. V. Vasconcelos, On the divisor class group of Rees algebras, 
{\em J. Algebra} {\bf 93} (1985), 182--188.} 







\bibitem{syl1}{J. Hong, A. Simis and W. V. Vasconcelos, The homology
of two-dimensional elimination, {\em J. Symbolic Computation}
{43} (2008),  275--292.}

\bibitem{syl2}{J. Hong, A. Simis and W. V. Vasconcelos, On the equations
of almost complete intersections,  {\em Bull. Braz. Math.  Soc.} {\bf 43} (2012), 171-199.}

\bibitem{Huck87}{S. Huckaba,  Reduction number for ideals of higher analytic spread, 
{\em Math. Proc. Camb. Phil. Soc.} {\bf 102} (1987), 49--57.} 


\bibitem{Huc96}{S. Huckaba, A $d$-dimensional extension of a lemma of
Huneke's and formulas for the Hilbert coefficients, {\em Proc. Amer.
Math. Soc.} {\bf 124} (1996), 1393-1401.}


\bibitem{Johnson}{M. Johnson, Depth of symmetric algebras of certain ideals, {\em Proc. 
Amer. Math. Soc.} {\bf 129} (2001), 1581--1585.}


\bibitem{KPU}{A. Kustin, C. Polini and B. Ulrich, Rational normal
scrolls and the defining equations of Rees algebras,
{\em J. Reine Angew. Math.} {\bf 650} (2011), 23--65.}

\bibitem{MPV12}{F. Mui\~nos and F. Planas-Vilanova, The equations of Rees algebras of equimultiple ideals
of deviation one, {\em Proc. Amer. Math. Soc.}, to appear.}


\bibitem{Northcott}{D. G. Northcott, A homological investigation of a certain residual ideal, 
{\em Math. Annalem} {\bf 150} (1963), 99--110.} 


\bibitem{Rossi83}{M. E. Rossi,  A note on symmetric algebras which are Gorenstein, {\em Comm. Algebra} 
{\bf 11} (1983), 2575--2591.}



\bibitem{Simis04}{A. Simis, Cremona transformations and some related algebras, {\em J. Algebra} {\bf 280} 
(2004), 162--179.}


\bibitem{Trung98}{N. V. Trung, The Castelnuovo regularity of the Rees algebra and associated graded ring,
{\em Trans. Amer. Math. Soc.} {\bf 350} (1998), 2813--2832.}








\bibitem{alt}{W. V. Vasconcelos, {\em Arithmetic of Blowup Algebras},
  London Math. Soc., Lecture Note Series {195}, Cambridge
  University Press, 1994.}


\bibitem{icbook}{W. V. Vasconcelos, {\em Integral Closure}, Springer Monographs in Mathematics, New York, 2005. }





\end{thebibliography}
\end{document}